\documentclass[A4paper,12pt]{article}
\usepackage{latexsym}
\usepackage{mathrsfs}
\usepackage{amssymb}
\usepackage{amscd}
\usepackage[dvips]{graphicx}                  
\newtheorem{definition}{Definition}
\newtheorem{theorem}{Theorem}

\newtheorem{lemma}{Lemma}

\newenvironment{proof}[1][Proof]{\textbf{#1.} }{\ \rule{0.5em}{0.5em}}
\date{}

\long\def\symbolfootnote[#1]#2{\begingroup%
	\def\thefootnote{$\;$}\footnote[#1]{$^*$#2}\endgroup}
\begin{document}
	
	\title{On some generalizations \\of the Halpern-L\"auchli theorem}
	\author{Joanna Jureczko}
	
	\maketitle
	
	\symbolfootnote[2]{Mathematics Subject Classification (2010): 54C30, 03E05, 03E40, 28A20.  
		
		\hspace{0.2cm}
		}
	
	\begin{abstract}
		In this paper, several generalizations of the classical Halpern-L\"{a}uchli Theorem are proven for Marczewski and Ellentuck structures using only combinatorial methods.
	\end{abstract}	
			
		\section{Introduction}
		
The Classical Halpern-Läuchli Theorem (\cite{HL}) concerns products of finitely many trees of height $\omega$ which have finite branches but do not have terminal nodes. More specifically, given a perfect tree $T$ and denoting the $n$-th level of $T$ by $T(n)$, the Halpern-Läuchli Theorem for $d<\omega$ (shortly $HL_{d}$) states that.
\\\\
\textbf{Theorem [Halpern-L\"auchli]}
\textit{If $(T_i \colon i < d)$ is a sequence of perfect trees, $A$ is an infinite subset of $\omega$ and 
$$\bigcup_{n \in A} \bigotimes_{i < d} T_i(n) = G_0 \cup G_1,$$
then there are $j< 2$,  an infinite subset $B$ of $A$  and downwards closed perfect subtrees $T'_i$ of $T_i, (i< \omega),$ with
$$\bigcup_{n \in B}\bigotimes_{i < d} T'_i(n) \subseteq G_j.$$}

Laver in \cite{RL} proved that the above result is true for $d = \omega$, which we will refer to as $HL_\omega$ in the subsequent parts of this paper. We can restate Laver's Theorem as follows.
\\\\
\textbf{Theorem 0.} 
\textit{If $f_i, (i < \omega)$ are continuous functions from the Hilbert cube $[0,1]^\omega$ into $[0,1]$, then there exist non-empty perfect sets $P_i \subseteq [0,1], (i < \omega)$ and $B \in [\omega]^\omega$ such that $(f_i\colon i \in B)$ is monotonic (and uniformly convergnet) on $\bigotimes_{i < \omega}P_i$.}
\\\\
The last theorem is also known as Harrington's one, (compare \cite[p. 481]{AS}).

The Halpern-L\"auchli theorem itself has received several generalizations, including versions for uncountable and measurable cardinals (see \cite{DH1, DH2}), as well as numerous applications in the proof of partition relation theorems (see the introduction in \cite{DH1, DH2}). The proofs of these results often rely on forcing methods, such as in \cite{TF, SS1}. For further discussions on the Halpern-L\"auchli theorem and its generalizations, refer to \cite{DH1, DH2, ST}.

In this paper, we will prove some generalizations of the Halpern-L"auchli Theorem in Marczewski and Ellentuck structures using only combinatorial methods. We will ensure clarity in our formulations of the results, basing them on the Halpern-L"auchli version given in Theorem 0.

It is worth emphasizing that the problem for "one dimension" was solved by S. Mazurkiewicz in 1920, but there is no bibliographical data supporting this fact.

Despite the existence of different versions of the Halpern-L\"auchli Theorem in the literature, no generalization of this theorem has been found for different "types" of measurability of functions $f_i$ given in Theorem 0.

It is generally known that there is a duality between Lebesgue measurability and the Baire property, see \cite{JO}. In fact, if we consider the set $X$ (which is not necessarily a topological space), among the subsets of $X$, we can consider the family $\mathcal{A}$ of collections (e.g., $\sigma-$field or Boolean algebra). Elements of this family are called \textit{large} sets. From the elements of $\mathcal{A}$, we choose the ones that make the ideal. Elements of this ideal are called \textit{small} sets.

In this paper, we provide considerations for three structures: the structure of Marczewski sets (Sacks forcing), the structure of Laver structure (Laver forcing), and the Ellentuck topology (Mathias forcing), about which we do not know if it makes a topology. For the Marczewski structure, $(s^0)-$sets will be equivalent to small sets and $(s)-$sets will be equivalent to large sets; for the Laver structure, $(l^0)-$sets will be equivalent to small sets and $(l)-$ sets will be equivalent to large sets. Finally, for the Ellentuck structure, $NR-$sets are the equivalents of small sets and $CR-$sets are the equivalents of large sets.

We have chosen these structures not accidentally because, on the one hand, they are very different, but on the other hand, they are very similar from a technical point of view. For all constructions considered in this paper, it is essential to use the Fusion Lemma (see \cite{TJ1, TJ, ST}), which was introduced in perfect sets forcing (i.e., Sacks forcing). It is worth adding that fusion in each mentioned structure runs in a different way.

Originally, we wanted to show our thesis only for Marczewski and Ellentuck structures, but the Laver structure is intuitively somewhere between these two structures. However, the reasonings presented below are restricted to Marczewski, Laver, and Ellentuck structures. It is also true for Silver-Prikry forcing (similar to Sacks forcing), Miller forcing (similar to Laver forcing), and Superset forcing, which is between Sacks and Laver forcing (see \cite{TJ1, JMS} for definitions of these notions).

Partial motivation for the considerations presented in this paper is based on \cite{AFP}. It has been noticed that measurability, $(s)-$measurability, and Ramsey-measurability are not equivalent concepts. It therefore seems essential to consider such generalizations of Theorem 0 in order to obtain further results in this direction.

Our earlier research on the Marczewski, Laver, and Ellentuck structures \cite{FJW} led us to solve the Kuratowski problem of 1935 \cite{KK1}. Both the results of \cite{FJW} and the methods of their proofs turned out to be useful in proving the generalizations of Halpern-L"auchli theorem proposed in this paper. The results presented in this paper concern the generalization of the Halpern-L\"auchli theorem (in the formulation of Theorem 0) for $(s)$- and $(l)$-measurability (where the concept of the Laver tree is replaced by a Laver-like tree - see constructions) and Ramsey-measurability.

As already mentioned in \cite{FJW}, the key tool used in the proofs of the theorems is the Fusion Lemma. Thus, we can assume that the results given here are true for all structures where the Fusion Lemma holds.

The proofs of the results where the Fusion Lemma is applied are not easy to write. For example, the proofs of Lemma 2 and Lemma 3 resemble the Sacks and Laver forcing, respectively (see \cite{TJ1}), while the proof of Lemma 5 resembles Mathias forcing (compare \cite{BJ}). However, in each case, the construction is different, so we decided to have quite extensive expanded Sections 2 containing definitions and previous results and detailed proofs of the main results in Section 3, which were divided into two parts, separately for tree structures and for the Ellentuck structure. The main results in Section 3 (i.e., Theorem 1 and Theorem 2) were preceded by auxiliary lemmas for easier and clearer study of the paper.

The structure of the paper is as follows: In Section 2, we give definitions and previous results needed in the further parts of this paper (some of them are rewritten from \cite{FJW}). For definitions and facts not cited here, we refer the reader to \cite{TJ, TJ1, RE}. In Section 3, divided into two subsections, we show auxiliary lemmas and the main results.
	
	\section{Definitions and previous results}
	
	\subsection{Tree ideals}
	
	Let $K\subseteq \omega$ be a set, (finite or infinite).
	A set $T\subseteq K^{<\omega}$ is called a \textit{tree} iff $t\upharpoonright n \in T$ for all $t \in T$ and $n\leqslant |t|$, (i.e. $T$ is closed downwards under initial segments). It is assumed that trees have no terminal nodes.
	
	Let $\mathbb{T}$ means a family of all tress. 
	For each $T \in \mathbb{T}$ and $t \in T$ the set 
	$$split(t, T) = |\{n \in K \colon t^\smallfrown n \in T\}|$$ denotes the number of successors of nodes in $T$.

	\begin{definition}
		A tree $T$ is called 
		\begin{enumerate}
			\item Sacks or perfect tree iff $K=\{0,1\}$ and for each $t \in T$ there is and $s \in T$ such that $t \subseteq s$ and $split(t, T) = 2$,
			\item Laver tree iff $K=\omega$ and there is $s \in T$ such that for each $t \in T$ 
			\begin{itemize}
			\item [(a)] either $t \subseteq s$ or $s \subseteq t$,	
			\item [(b)] $split(t, T)$ is infinite for each $t \in T$.  
			\end{itemize}			
		\end{enumerate}
	\end{definition}

\noindent
We denote by $\mathbb{S}$ ($\mathbb{L}$) the family of all Sacks (Laver) trees, respectively.

In the results below, we use a Laver-like structure of a tree (which is rather close to superset or Miller trees) because, for each $n \in \omega$, we divide the $n$-th level of nodes into $(n + 1)$ subsets of nodes.
\\

We say that $t\in Lev_n (T)$, (i.e. $t$ belongs to $n-$level of $T$) iff there are $n-$splits below $t$.
	
	Let $$[T] = \{x \in K^\omega \colon \forall_{n \in \omega}\ x \upharpoonright n \in T \}$$
	be the set of all infinite paths through $T$.
	\\
	Notice that $[T]$ is closed in the Baire space $K^\omega$, (see e.g. \cite{TJ}).
	
	By $stem(T)$ we mean a node $t \in T$ such that $split(s, T) = 1$ and $split(t, T) >~1$ for any $s \varsubsetneq t$.
	\\
	
	The ordering on $\mathbb{S}$ is defined as follows
	$Q\leqslant T $ iff $ Q \subseteq T$
	and 
	$$Q \leqslant_n T \textrm{ iff } Q \leqslant T \textrm{ and  any node of $n-$level of $T$ is a node of $n-$level of $Q$}.$$
	\indent
	If $T \in \mathbb{L}$, then $\{s \in [T] \colon stem(T) \subseteq s \}$, (i.e. the part of $T$ above the $stem(T)$),  can be enumerated as follows:
	$$s^T_0 = stem(T), s^T_1, ..., s^T_n, ...\ .$$
	Thus we can define the ordering on $\mathbb{L}$ in the following way: let $Q, T \in \mathbb{L}$, $Q \leqslant T$ iff $Q \subseteq T$
	and
	$$Q \leqslant_n T \textrm{ iff } stem(Q) = stem(T) \textrm{ and } s^Q_i = s^T_i \textrm{ for all } i = 0, ..., n.$$

	We say that a set $A \subseteq K^\omega$ is a \textit{$(t)-$set} iff 
	$$\forall_{T \in \mathbb{T}}\ \exists_{Q \in \mathbb{T}}\ Q \subseteq T \wedge ([Q] \subseteq A \vee [Q]\cap A =\emptyset). $$
	We say that a set $A \subseteq K^\omega$ is a \textit{$(t^0)-$set} iff 
	$$\forall_{T \in \mathbb{T}}\ \exists_{Q \in \mathbb{T}}\ Q \subseteq T \wedge  [Q]\cap A =\emptyset. $$
	
	Throughout the paper we assume that a set of trees is a $(t^0)-$set iff its set of infinite paths is a $(t^0)-$set in $K^\omega$.
	Thus, we will denote by $\mathbb{S}^0$ $(\mathbb{L}^0)$ the family of all $(s^0)-$sets ($(l^0)-$sets), respectively. 
	
	The fact that $(s^0)$- and $(l^0)$-sets are $\sigma$-ideals in $2^\omega$ and $\omega^\omega$, respectively, is an application of the Fusion Lemma (see Section 2.4 below).
	
	
	For further consideration, unless otherwise stated, $\mathbb{T}$ and $\mathbb{T}^0$ refer to the $\sigma-$ideals of Sacks trees, i.e., the $(s)-$ and $(s_0)-$trees, respectively, and Laver trees, i.e., the $(l)-$ and $(l^0)-$trees, respectively. Then $(t), (t^0)$, and $K$ will be determined accordingly to these structures.

	\subsection{Ellentuck topology}
	
	The Ellentuck topology $[\omega]^{\omega}_{EL}$ on $[\omega]^\omega$ is generated by sets of the form
	$$[a, A] = \{B \in [A]^\omega \colon a \subset B \subseteq a \cup A\},$$
	where $a \in [\omega]^{<\omega}$ and $A \in [\omega]^\omega$. We call such sets \textit{Ellentuck sets, (shortly $EL-$sets).} Obviously $$[a, A] \subseteq [b, B] \textrm{ iff } b \subseteq a \textrm{ and } A \subseteq B.$$ 
	
	A set $M \subseteq [\omega]^\omega$ is \textit{completely Ramsey}, (shortly \textit{$CR-$set}), if for every $[a, A]$ there exists $B \in [A]^\omega$ such that $[a, B] \subseteq M$ or $[a, B] \cap M = \emptyset.$ 
	A set $M \subseteq [\omega]^\omega$ is \textit{Ramsey null}, (shortly \textit{$NR-$set}), if for every $[a, A]$ there exists $B \in [A]^\omega$ such that $[a, B] \cap M = \emptyset.$

	Notice that the application of the Fusion Lemma (see Section 2.4 below) implies that all $NR$-sets form a $\sigma$-ideal in $[\omega]^{\omega}_{EL}$, which we denote by $\mathbb{NR}$.
\\

	\noindent The following fact will be used in further considerations for simplifying the notation. 
	\\\\
	\noindent
	\textbf{Fact 1 (\cite{JB})} \textit{Let $M$ be an open and dense set (in the sense of Ellentuck topology). Then for each $A \subseteq [\omega]^\omega$ and for each $a \in [\omega]^{<\omega}$ there exists $B \subseteq [\omega]^\omega$ such that $B \subseteq A$ the set $[\emptyset, B \cup a] \subseteq M$.} 
	
	\subsection{Kuratowski partitions}
	
		Let $A \subseteq K^\omega$.
		A partition $\mathcal{F}$ of $A$  into $(t^0)-$sets is called \textit{Kuratowski partition} if $\bigcup\mathcal{F}'$ is a $(t)-$set for any $\mathcal{F}' \subseteq \mathcal{F}$.
	\\
	
	Let $M \subseteq [\omega]^\omega$.
	A partition $\mathcal{F}$ of $M$  into $NR-$sets is called  \textit{Kuratowski partition} if $\bigcup\mathcal{F}'$ is a $CR-$set for any $\mathcal{F}' \subseteq \mathcal{F}$.
\\

\noindent
\textbf{Fact 2 (\cite{FJW})}
\textit{\begin{enumerate}
		\item No $A \in P(K^\omega)\setminus \mathbb{T}^0$ admits  Kuratowski partition.
		\item No open and dense set $M\in  P([\omega]^{\omega})\setminus \mathbb{NR}$ admits  Kuratowski partition.
\end{enumerate}}

\noindent
\textbf{Fact 3 (\cite{FJW})}
\textit{Let $A \in P(2^\omega) \setminus \mathbb{S}^0$. For any  partition $\mathcal{F}$ of $A$ into $(s^0)-$sets and for any perfect tree $T\in \mathbb{S}$ there exists a perfect subtree $Q \leqslant T$  such that the family
	$$\mathcal{F}_{[Q]} = \{F \cap [Q] \colon F \in \mathcal{F}\}$$ has cardinality continuum.}
\\
\\
\textbf{Fact 4 (\cite{FJW})}
\textit{Let $A \in P(\omega^\omega) \setminus \mathbb{L}^0$. For any  partition $\mathcal{F}$ of $A$ into $(l^0)-$sets and for any Laver tree $T\in \mathbb{L}$ there exists a Laver subtree $Q \leqslant T$  such that the family
	$$\mathcal{F}_{[Q]} = \{F \cap [Q] \colon F \in \mathcal{F}\}$$ has cardinality continuum.}
\\\\
\textbf{Fact 5 (\cite{FJW})}
\textit{Let $M \in P([\omega]^\omega)\setminus \mathbb{NR}$ be an open and dense set (in the sense of Ellentuck topology). For any partition $\mathcal{F}$ of $M$ into $NR-$sets and for any $[a, A] \subseteq [\omega]^{\omega}_{EL} $ there exists $[b, B] \subseteq [a, A]$ such that the family $$\mathcal{F}_{[b, B]} = \{F \cap [b, B] \colon F \in \mathcal{F}\}$$ has cardinality continuum.}
\\

A family  $\{T_i \colon i \in I\}$ of subtrees of $\mathbb{T}$ is called \textit{$(t)-$additive} if for any $J \subseteq I$ the set $\bigcup\{T_i \colon i \in J\}$ is a $(t)-$set.
\\
\\
	\textbf{Fact 6 (\cite{AFP})}
If $\{S_i \colon i \in I\}$  is a disjoint $(s)-$additive family contained in $\mathbb{S}^0$, then $\bigcup\{S_i \colon i \in I\}$  belongs to $\mathbb{S}^0$.
\\

A family  $\{M_i \colon i \in I\}$ of sets $M_i \subseteq [\omega]^\omega$  is called \textit{$CR-$additive} if for any $J \subseteq I$ the set $\bigcup\{M_i \colon i \in J\}$ is a $CR-$set.
\\\\
	\textbf{Fact 7 (\cite{AFP})}
If $\{M_i \colon i \in I\}$ of sets $M_i \subseteq [\omega]^\omega$ is a disjoint $CR-$additive family contained in $\mathbb{NR}$, then $\bigcup\{M_i \colon i \in I\}$  belongs to $\mathbb{NR}$.

	\subsection{Fusion Lemma}
	Let $\mathbb{T}$  be the family of all trees.  
	A sequence $\{T_n\}_{n \in \omega}$ of trees such that 
	$$T_0 \geqslant_0 T_1 \geqslant_1 ... \geqslant_{n-1} T_n \geqslant_n ...$$
	is called a \textit{fusion sequence}.
	\\
	\\
	\textbf{Fact 8 (\cite{TJ1})} If $\{T_n\}_{n\in \omega}$ is a fusion sequence then $T = \bigcap_{n\in \omega}T_n$, (the fusion of $\{T_n\}_{n \in \omega}$), belongs to $\mathbb{T}$.
	\\ 
	\\
	A sequence $\{[a_n, A_n]\}_{n \in \omega}$ of $EL-$sets  is called a \textit{fusion sequence} if it is infinite and
	\\(1) $\{a_n\}_{n \in \omega}$ is a nondecreasing sequence of integers converging to infinity,
	\\(2) $A_{n+1} \in [a_n, A_n]$ for all $n \in \omega$.
	\\
	\\
	\textbf{Fact 9 (\cite{TJ1})} If $\{[a_n, A_n]\}_{n \in \omega}$ is a fusion sequence then $$[a, A] = \bigcap_{n\in \omega}[a_, A_n] =  [a,\bigcap_{n\in \omega}  A_n],$$ (the fusion of $\{[a_n, A_n]\}_{n \in \omega})$, is an $EL-$set.

	\subsection{$(t)-$  and $CR-$measurable functions}
	
	Let $K$ be as in Section 2.1.
	A function $f \colon K^\omega \to \mathbb{R}$ is called \textit{$(t)-$measurable} if for any open set $U \subseteq \mathbb{R}$ the set $f^{-1}(U)$ is a $(t)-$set.
	\\\\	
	\textbf{Fact 10 (\cite{AFP})} \textit{For any metric space $X$ and for any function $f\colon \mathbb{R} \to X$, f is $(s)-$measurable iff for every perfect set $P$ there exists a perfect set $Q \subseteq P$, such that $f|_Q$ is continuous.} 
	\\
	
	A function $f \colon[\omega]^\omega \to \mathbb{R}$, is called $CR-$measurable if for any open set $U \subseteq \mathbb{R}$, the set $f^{-1}(U)$ is a $CR-$set.
	\\\\
	\textbf{Fact 11 (\cite{AFP})} \textit{For any metric space $X$ and for any function $f\colon [\omega]^\omega \to X$, f is $CR-$measurable iff for every $[a, A]$ there exists infinite $B \subseteq A$, such that $f|_{[a, B]}$ is continuous, ($[\omega]^\omega \subseteq \{0,1\}^\omega$ and continuity is regarded in the subspace topology).}

	\section{Main results}
	
	\subsection{$HL_{\omega}$ in tree structures}
	
	We will begin this subsection with auxiliary lemmas. Lemma 1 and Lemma 2 are generalizations of Fact 6 and Fact 10, respectively.
	
	\begin{lemma}
		If $\{T_i \colon i \in I\}$ is a disjoint $(t)-$additive family contained in $\mathbb{T}^0$ then $\bigcup\{T_i \colon i \in I\}$ belongs to $\mathbb{T}^0.$ 
	\end{lemma}
	
	\begin{proof} Suppose that $\bigcup\{T_i \colon i \in I\}$ does not belong to $\mathbb{T}_0$.
		Let $T \subseteq \bigcup\{T_i \colon i \in I\}$ be  a $(t)-$tree. Let $\{Q_\alpha \colon \alpha \in m^\omega\}$ be a set of all $(t)-$subtrees of $T$, where $m=2$ for $(s)-$trees and $m=\omega$ for $(l)-$trees.  By Fact 3 and Fact 4, respectively, there exist distinct $i_\alpha, j_\alpha \in \{i_\beta, j_\beta \in I \colon \beta < \alpha\}$ such that 
		$$[Q_\alpha] \cap[T_{i_\alpha}] \not = \emptyset \textrm{ and }
		[Q_\alpha] \cap[T_{j_\alpha}] \not = \emptyset.$$ 
		Let $J = \{j_\alpha \colon \alpha \in m^\omega\}$. Then for all $\alpha \in m^\omega$ $$[Q_\alpha] \cap \bigcup\{T_i \colon i \in J\}   \not = \emptyset $$ and $$ [Q_\alpha] \cap \bigcup\{T_i \colon i \in I \setminus J\}  \not = \emptyset.$$
		Hence $\bigcup\{T_i \colon i \in J\}$ is not a $(t)-$set. A contradiction.
	\end{proof}
	\\

As mentioned in Section 2.1, for the constructions given in Lemma 2 and Lemma 3, we use a Laver-like tree structure. On the $n$-th level of nodes, we divide them into $(n+1)$ subsets of nodes for each $n \in \omega$.
	
	\begin{lemma}
		Let $X$ be a metric space and let $f \colon K^\omega \to X$ be a function. Then, $f$ is $(t)-$measurable iff for any $(t)-$tree $T \in \mathbb{T}$ there exists a $(t)-$subtree $Q \subseteq T$ such that $f|_{[Q]}$ is continuous. 
	\end{lemma}
	
	\begin{proof}
		Assume that $f$ is $(t)-$measurable and $T\in \mathbb{T}$ is a $(t)-$tree.
		We will construct by induction $(t)-$trees $\{T_h \colon h \in m^n, m, n \in \omega\}$. The first and inductive steps are the same. 
		Assume that we have chosen $T_h, (h \in m^{n-1})$ such that $diam(f(T_h)) \leqslant \frac{1}{m^n}$. Consider an open cover of $X$
		$$\mathcal{U}_n = \{V \subset X \colon diam (V) < (\frac{1}{2})^n\}.$$
		By the Stone Theorem, (see \cite{RE}), the family $\mathcal{U}_n$ has a $\sigma-$discrete refinement $\bar{\mathcal{U}}_n$.
		Since $f$ is $(t)-$measurable the family $\{f^{-1}(U) \colon U \in \bar{\mathcal{U}}_n\}$ is disjoint and $(t)-$additive.
		By Lemma 1, $\bigcup\{f^{-1}(U) \colon U \in \bar{\mathcal{U}}_n\}$ belongs to $\mathbb{T}^0$. Hence, for any $h \in m^{n-1}$ we can choose 
		$U \in \bar{\mathcal{U}}_n$ such that $f^{-1}(U)\cap [T_h] \in \mathbb{T}\setminus \mathbb{T}^0$. in this set we can construct 
	 $(t)-$trees fulfilling the properties: for all distinct $h, h' \in m^n$, $m,n \in \omega$ 
		\begin{itemize}
			\item [(1)] $T_h \leqslant_n T_{h|(n-1)}$, i.e. $[T_h] \subseteq [T_{h|(n-1)}]$;
			\item [(2)] $[T_h]\cap [T_{h'}] = \emptyset$;
			\item [(3)]  $diam(f(T_h)) \leqslant \frac{1}{m^n}$.
		\end{itemize}
		Since the constructions of $(s)-$trees and  $(l)-$trees run in the different ways, we construct them separately. 
		\\
		
		Case 1. ($m=2$, the construction of $(s)-$trees)
		\\
		Assume that for some $n \in \omega$ we have constructed the family $\{T_h \colon h \in 2^n\}$ of properties $(1)-(3)$.
		We will construct $(n+1)-$step.
		
		Fix $h \in 2^n$. Let $B_h(n)$ be a set of all nodes of $n-$level of $T_h$.
		Since $split(t, T) = 2$, we can divide $B_h(n)$ into disjoint sets $B_{h^\smallfrown 0}(n+1)$ and $ B_{h^\smallfrown 1}(n+1)$ which will be $(n+1)-$levels of $(s)-$trees  $T_{h^\smallfrown 0}, T_{h^\smallfrown 1} \leqslant_n T_h$, respectively.
		
		Define
		$$t\in B_{h^\smallfrown 0}(n+1) \textrm{ iff } \exists_{x \in [T_h]}\   t \in x  \textrm{ and } x(n) =0,$$
		(where $x(n)$ is the $n-$th element of $x \in 2^\omega$ considered as $x \colon \omega \to \{0, 1\}$).
		In  the similar way we define  $B_{h^\smallfrown 1}(n+1)$.
		
		Now, choose $(s)-$subtrees $T_{h^\smallfrown 0}, T_{h^\smallfrown 1} \leqslant _n T_h$ such that
		$$\{x \in [T_h] \colon \exists_{t \in B_{h^\smallfrown 0}(n+1)\setminus B_{h^\smallfrown 1}(n+1)}\  t \in x\} \subseteq [T_{h^\smallfrown 0}] $$
		$$\{x \in [T_h] \colon \exists_{t \in B_{h^\smallfrown 1}(n+1)\setminus B_{h^\smallfrown 0}(n+1)}\  t \in x\} \subseteq [T_{h^\smallfrown 1}]$$
		$$diam(f(T_{h^\smallfrown 0})) \leqslant \frac{1}{2^{n+1}} \textrm{ and } diam(f(T_{h^\smallfrown 1})) \leqslant \frac{1}{2^{n+1}}.$$
		The construction of $(n+1)-$step for $(s)-$trees is complete.
		
		Case 2. ($m=n$, the construction of $(l)-$trees)
		\\
		Assume that for some $n \in \omega$ we have constructed the family $\{T_h \colon h \in n^n\}$ of properties $(1)-(3)$.
		We will construct now $(n+1)-$step.
		
		Fix $h \in n^n$. Let $B_h(n)$ be a set of all nodes of $n-$level of $T_h$.
		Since $split(t, T) = \omega$, we can divide $B_h(n)$ into disjoint  $k$ sets $B_{h^\smallfrown k}(n+1)$ for $k=0, 1, ..., n$ which will be $(n+1)-$levels of $(l)-$trees  $T_{h^\smallfrown k} \leqslant_n T_h$, respectively.
		
		Define
		$$t\in B_{h^\smallfrown 0}(n+1) \textrm{ iff } \exists_{x \in [T_h]}\   t \in x  \textrm{ and } x(n) =0,$$
		$$t\in B_{h^\smallfrown k}(n+1) \textrm{ iff } \exists_{x \in [T_h]}\   t \in x  \textrm{ and } x(n) =k, t \not \in \{B_{h^\smallfrown l}(n+1)\colon l<k\}$$
		(where $x(n)$ is the $n-$th element of $x \in \omega^\omega$).

		Now, choose $(l)-$subtrees $T_{h^\smallfrown k} \leqslant _n T_h$ for $k =0, 1, ..., n$ such that
		$$\{x \in [T_h] \colon \exists_{t \in B_{h^\smallfrown k}(n+1)}\  t \in x\} \subseteq [T_{h^\smallfrown 0}] $$
		and
		$$diam(f(T_{h^\smallfrown k})) \leqslant \frac{1}{n+1^{n+1}}.$$
		The construction of $(n+1)-$step for $(l)-$trees is complete.
		
		Now, take $Q = \bigcap_{n < \omega} \bigcup_{h \in m^n}T_h$. By Fact 8, $Q$ is a $(t)-$tree. It is easy to see that $f|_{[Q]}$ is continuous.
		
		Assume now, that for any $(t)-$tree $T \in \mathbb{T}$  there is a  $Q\subseteq T$ such that $f|_{[Q]}$ is continuous.
		
		Let $U \subseteq X$ be an open set. Then, $f^{-1}(U)\cap [Q]$ is open in $[Q]$ and hence $f^{-1}(U)\cap [Q]$ or its complement in $[Q]$  contains a $(t)-$subtree. This completes the proof. 
	\end{proof}
	\\

The next lemma is a "one-dimensional" version of the main result in this subsection, given in Theorem 1.

	\begin{lemma}
		Let $f_n \colon K^\omega \to \{0,1\}$ be a sequence of $(t)-$measurable functions, $n \in \omega$. Then, there exists a $(t)-$tree $Q \subseteq K^{<\omega}$ and a subsequence $\{f_{n_k}\}$ which is uniformly convergent on $[Q]$. 
	\end{lemma}
	
	\begin{proof}
		Using Lemma 2, it is enough to assume that the functions $f_n, (n \in \omega)$ are continuous on some $(t)-$tree $T$.
		
		We will construct a family of $(t)-$trees $\{T^{\varepsilon_{n_k}}_h \colon h \in m^n, \varepsilon_{n_k} \in \{0, 1\}\}$ of $T$ and a subsequence $\{f^{\varepsilon_{n_k}}_{n_k}\}$ of $\{f_n\}$ of the following properties: for all distinct $h, h' \in m^k$	
		\begin{itemize}
			\item [(1)] $T^{\varepsilon_{n_k}}_h \leqslant_k T^{\varepsilon_{n_{k-1}}}_{h|(k-1)}$, i.e. $[T^{\varepsilon_{n_k}}_h] \subseteq [T^{\varepsilon_{n_{k-1}}}_{h|(k-1)}]$;
			\item [(2)] $[T^{\varepsilon_{n_k}}_h]\cap [T^{\varepsilon_{n_k}}_{h'}] = \emptyset$;
			\item [(3)] $T^{\varepsilon_{n_k}}_h\leqslant_k (f^{\varepsilon_{n_k}}_{n_k})^{-1}(\{\varepsilon_{n_k}\})$.
		\end{itemize}
	
	The notation $\varepsilon_{n_k}$ in conditions $(1)-(3)$ is essential because $\varepsilon_{n_k}$ can be differnt at each level, i.e. $(\varepsilon_{n_k})_{n_k \in \omega}$ form a $0-1$ sequences.
		
		Since the constructions of $(s)-$trees and  $(l)-$trees run in the different ways, we construct them separately. 
		\\	
		
		Case 1. (m=2, the construction of $(s)-$trees)
		\\
		Let $f^{\varepsilon_{n_0}}_{n_0}$ be an arbitrary element of $\{f_n\}$. The first and inductive steps are the same. Assume that for some $k\in \omega$ there is constructed the subsequence $\{f^{\varepsilon_{n_k}}_{n_k}\}$ and the family of $(s)-$subtrees $$\{T^{\varepsilon_{n_k}}_h \colon h \in 2^k, \varepsilon_{n_k} \in \{0,1\}\}$$ fulfilling $(1)-(3)$. (At least one of  $T^{\varepsilon_{n_k}}_{h}$ is nonempty, so one can continue the construction only for such $(s)-$trees).
		
		Fix $h\in 2^k$ for which $T^{\varepsilon_{n_k}}_{h}$ is nonempty. Choose a function $f^{\varepsilon_{n_{k+1}}}_{n_{k+1}}$ from $\{f_n\}$ such that $$[T^{\varepsilon_{n_k}}_h] \cap (f^{\varepsilon_{n_{k+1}}}_{n_{k+1}})^{-1}(\{\varepsilon_{n_k}\}) \not = \emptyset.$$
		Let $B^{\varepsilon_{n_k}}_h(k)$ be a set of all nodes of $k-$level of $T^{\varepsilon_{n_k}}_h$.
		Since $split(t, T) = 2$, we can divide $B^{\varepsilon_{n_k}}_h(k)$ into disjoint sets $B^{\varepsilon_{n_{k+1}}}_{h^\smallfrown 0}(k+1)$ and $ B^{\varepsilon_{n_{k+1}}}_{h^\smallfrown 1}(k+1)$ which will be $(k+1)-$levels of $(s)-$trees  $T^{\varepsilon_{n_{k+1}}}_{h^\smallfrown 0}, T^{\varepsilon_{n_{k+1}}}_{h^\smallfrown 1} \leqslant_k T^{\varepsilon_{n_k}}_h$, respectively.
		
		Define
		$$t\in B^{\varepsilon_{n_{k+1}}}_{h^\smallfrown 0}(k+1) \textrm{ iff } \exists_{x \in [T_h]}\   t \in x  \textrm{ and } x(k) =0,$$
		(where $x(k)$ is the $k-$th element of $x \in 2^\omega$ considered as $x \colon \omega \to \{0, 1\}$).
		In the similar way we define  $B^{\varepsilon_{n_{k+1}}}_{h^\smallfrown 1}(k+1)$.
		
		Now, choose $(s)-$subtrees $T^{\varepsilon_{n_{k+1}}}_{h^\smallfrown 0}, T^{\varepsilon_{n_{k+1}}}_{h^\smallfrown 1} \leqslant _k T^{\varepsilon_{n_{k}}}_h$ such that
		$$\{x \in [T^{\varepsilon_{n_{k}}}_h] \colon \exists_{t \in B^{\varepsilon_{n_{k+1}}}_{h^\smallfrown 0}(k+1)\setminus B^{\varepsilon_{n_{k+1}}}_{h^\smallfrown 1}(k+1)}\  t \in x\} \subseteq [T^{\varepsilon_{n_{k+1}}}_{h^\smallfrown 0}] $$
		$$\{x \in [T^{\varepsilon_{n_{k}}}_h] \colon \exists_{t \in B^{\varepsilon_{n_{k+1}}}_{h^\smallfrown 1}(k+1)\setminus B^{\varepsilon_{n_{k+1}}}_{h^\smallfrown 0}(k+1)}\  t \in x\} \subseteq [T^{\varepsilon_{n_{k+1}}}_{h^\smallfrown 1}]$$
		$$T^{\varepsilon_{n_{k+1}}}_{h^\smallfrown 0}\leqslant_{k+1} (f^{\varepsilon_{n_{k+1}}}_{n_{k+1}})^{-1}(\{\varepsilon_{n_{k+1}}\})$$
		and 
		$$T^{\varepsilon_{n_{k+1}}}_{h^\smallfrown 1}\leqslant_{k+1} (f^{\varepsilon_{n_{k+1}}}_{n_{k+1}})^{-1}(\{\varepsilon_{n_{k+1}}\}).$$
		The construction of $(k+1)-$step for $(s)-$trees is complete.

		Case 2. ($m=k$, the construction of $(l)-$trees)
		\\
		Let $f^{\varepsilon_{\varepsilon_{n_0}}}_{n_0}$ be an arbitrary element of $\{f_n\}$. The first and inductive steps are the same. Assume that for some $k\in \omega$ there is constructed the subsequence $\{f^{\varepsilon_{n_k}}_{n_k}\}$ and the family of $(l)-$subtrees $$\{T^{\varepsilon_{n_k}}_h \colon h \in 2^k, \varepsilon_{n_k} \in \{0,1\}\}$$ fulfilling $(1)-(3)$. (At least one of  $T^{\varepsilon_{n_k}}_{h}$ is nonempty, so one can continue the construction only for such $(l)-$trees).
		
		Fix $h\in k^k$ for which $T^{\varepsilon_{n_k}}_{h}$ is nonempty.
		Choose a function $f^{\varepsilon_{n_{k+1}}}_{n_{k+1}}$ from $\{f_n\}$ such that $$[T^{\varepsilon_{n_k}}_h] \cap (f^{\varepsilon_{n_{k+1}}}_{n_{k+1}})^{-1}(\{\varepsilon_{n_k}\}) \not = \emptyset.$$
		Let $B^{\varepsilon_{n_k}}_h(k)$ be a set of all nodes of $k-$level of $T^{\varepsilon_{n_k}}_h$.
		Since $split(t, T) = \omega$, we can divide $B^{\varepsilon_{n_k}}_h(k)$ into disjoint  $k$ sets $B^{\varepsilon_{n_{k+1}}}_{h^\smallfrown m}(k+1)$ for $m=0, 1, ..., k$ and  which will be $(k+1)-$levels of $(l)-$trees $T^{\varepsilon_{n_{k+1}}}_{h^\smallfrown m} \leqslant_k T_h$, respectively.
		
		Define
		$$t\in B^{\varepsilon_{n_{k+1}}}_{h^\smallfrown 0}(k+1) \textrm{ iff } \exists_{x \in [T^{\varepsilon_{n_k}}_h]}\   t \in x  \textrm{ and } x(k) =0,$$
		$$t\in B^{\varepsilon_{n_{k+1}}}_{h^\smallfrown k}(k+1) \textrm{ iff } \exists_{x \in [T^{\varepsilon_{n_k}}_h]}\   t \in x  \textrm{ and } x(k) =k, t \not \in \{B^{\varepsilon_{n_{k+1}}}_{h^\smallfrown l}(k+1)\colon l<k\}$$
		(where $x(n)$ is the $n-$th element of $x \in \omega^\omega$).

		Now, choose $(l)-$subtrees $T^{\varepsilon_{n_k}}_{h^\smallfrown m} \leqslant _n T^{\varepsilon_{n_k}}_h$ for $m =0, 1, ..., k$ such that
		$$\{x \in [T^{\varepsilon_{n_k}}_h] \colon \exists_{t \in B^{\varepsilon_{n_{k+1}}}_{h^\smallfrown k}(k+1)}\  t \in x\} \subseteq [T^{\varepsilon_{n_{k+1}}}_{h^\smallfrown 0}] $$
		$$\{x \in [T^{\varepsilon_{n_k}}_h] \colon \exists_{t \in B^{\varepsilon_{n_{k+1}}}_{h^\smallfrown k}(k+1)}\  t \in x\} \subseteq [T^{\varepsilon_{n_{k+1}}}_{h^\smallfrown 1}] $$
		$$T^{\varepsilon_{n_{k+1}}}_{h^\smallfrown 0}\leqslant_{k+1} (f^{\varepsilon_{n_{k+1}}}_{n_{k+1}})^{-1}(\{\varepsilon_{n_{k+1}}\})$$
		and
		$$T^{\varepsilon_{n_{k+1}}}_{h^\smallfrown 1}\leqslant_{k+1} (f^{\varepsilon_{n_{k+1}}}_{n_{k+1}})^{-1}(\{\varepsilon_{n_{k+1}}\}).$$
		The construction of $(k+1)-$step for $(l)-$trees is complete.
		
		By Fact 8,  $$Q^{(\varepsilon_{n_k})_{n_k \in \omega}} = \bigcap_{k \in \omega} \bigcup_{h \in m^k} T^{\varepsilon_{n_k}}_h$$ are $(t)-$trees.
		Let $$Q^\varepsilon = \bigcup \{Q^{(\varepsilon_{n_k})_{n_k \in \omega}}\colon  (\varepsilon_{n_k})_{n_k \in \omega} \to \varepsilon \},$$
		where $\varepsilon \in \{0, 1\}$.
		One of the sets $Q^0$ or $Q^1$ is nonempty. Then, $$\{f^{\varepsilon_{n_k}}_{n_k}\colon \varepsilon_{n_k} = 0\}$$ is uniformly convergent on $[Q^0]$ or $$\{f^{\varepsilon_{n_k}}_{n_k}\colon \varepsilon_{n_k} = 1\}$$ is uniformly convergent on $[Q^1]$.
	\end{proof}
	
	\begin{theorem}
		Let $f_n \colon K^\omega \to [0,1]$ be a sequence of $(t)-$measurable functions. Then there are $(t)-$trees $\{T_k \colon k \in \omega\}$ and a subsequence $$\{f_{n_k} \colon k \in \omega\}$$ which is uniformly convergent on $\bigotimes_{k \in \omega} [T_k]$.
	\end{theorem}
	
	\begin{proof}
		Since $f_n \colon K^\omega \to [0,1]$ for all $n \in \omega$, are $(t)-$measurable, there exist $(t)-$measurable functions $g_n \colon K^\omega \to \{0,1\}^\omega$ such that $f_n = g_n \circ \varphi$, where $\varphi \colon \{0,1\}^\omega \to [0,1]$ is the canonical continuous "onto" mapping (consider elements of $[0,1]$ in their binary expansions). 
		Thus, for completing the proof it is enough to find $(t)-$trees  $T_k$ for $ k \in \omega$ and  a subsequence $\{g_{n_k}\}$ which is uniformly convergent on $\bigotimes_{k< \omega} T_k$.
		Since $(\{0,1\}^\omega)^{K^\omega}$ can be "identified" by an exponential law with $(\{0,1\}^{K^\omega})^\omega$ 
		 one can consider each $g_n$ to be in the form $g_n = \{ h^k_n \colon k\in \omega\}$, where $h^k_n \colon K^\omega \to \{0,1\}$ are $(t)-$measurable.
		 By Lemma 2, for each $h^k_n$ and for each $(t)-$tree $T$ there is $T'\leqslant T$ such that $[T']\subseteq [T]$ and $h^k_n|_{[T']}$ is continuous. 
		
		Now, a family ${{g_{m,i}} \colon m \in M_i, i< \omega}$ will be constructed (for some set of indices $M_i\subseteq \omega$ of the subsequences $\widehat{h}^i_{m}$ chosen at each step in the construction) of infinite subsequences of ${g_n}$ in the following way.

		Using Lemma 3 on the sequence ${h^0_n}$ (i.e., the first coordinates of ${g_n}$), we obtain a $(t)$-tree $T^{\varepsilon_0}0$ and a subsequence ${\widehat{h}^0{m}}$ that uniformly converges on $T^{\varepsilon_0}0$, where $\varepsilon_0 \in {0,1}$. Then we obtain $$g{m,0} = ( \widehat{h}^0_{m}|{T^{\varepsilon_0}0}, h^1_m, ..., h^k_m,... ).$$ Now we use Lemma 3 on the sequence ${h^1_m}$ (the second coordinates of $g{m,0}$), and similarly as above, we obtain a $(t)$-tree $T^{\varepsilon_1}1$ and a subsequence ${\widehat{h}^1{m}}$ that uniformly converges on $T^{\varepsilon_1}1$, where $\varepsilon_1 \in {0,1}$. Thus, we can define the subsequence ${g{m,1}}$ of ${g{m,0}}$, which is of the form $$g_{m,1} = ( \widehat{h}^0_{m}|_{T^{\varepsilon_0}0}, \widehat{h}^1_m|{T^{\varepsilon_1}1}, h^2_m, ..., h^k_m,... ).$$ Using $\omega$ times Lemma 3 for each ${h^k_m}$, we can obtain the family $${{g{m,i}} \colon m \in M_i, i< \omega}$$ of infinite subsequences with the following properties:
		\begin{itemize}
			\item[(a)] $\{g_{m, i}\}$ is a subsequence of $\{g_{m, i-1}\}$, i.e. $M_i \subseteq M_{i-1}$ for all $i \in\omega$,
			\item[(b)] $g_{m,i} = \langle \bar{h}^k_m \colon k< \omega\rangle$ such that
			$$\bar{h}^k_m = \left\{ \begin{array}{ll}
			\widehat{h}^l_{m}|_{T^{\varepsilon_l}_l} & \textrm{for $l\leqslant i$}\\
			h^l_m & \textrm{for $l>i$}\\
			\end{array} \right.$$
			where $\{\widehat{h}^l_{m}\}$ is a subsequence of $\{\widehat{h}^{l-1}_{m}\}$ uniformly convergent on $T^{\varepsilon_l}_l, \varepsilon_l \in \{0,1\}$ and $h^l_m$ is an l-th coordinate of $g_{m, l}$.
		\end{itemize}	
		Now, from each subsequence $\{g_{m,i}\}$ one chooses its $i-$th element and define the sequence $\{g_{i,i} \colon i \in \omega\}$. Now, divide $\{g_{i,i}\}$ into two sequences (one of them may not exist):
		$\{g^0_{i,i}\}$ iff $\{\bar{h}^i_i\}$ are convergent uniformly on $ T^{\varepsilon_i}_{i}$ for $ \varepsilon_i = 0$ and 	$\{g^1_{i,i}\}$ iff $\{\bar{h}^i_i\}$ are convergent uniformly on  $ T^{\varepsilon_i}_{i}$ for $ \varepsilon_i = 0$. Then $\{g^0_{i,i}\}$ is uniformly convergent on $\bigotimes_{i< \omega} [T^0_{i}]$ or $\{g^1_{i,i}\}$ is uniformly convergent on $\bigotimes_{i< \omega} [T^1_{i}]$. 
	\end{proof}
	
	\subsection{$HL_\omega$ in Ellentuck structures}

	This subsection incorporates the versions of Lemma 1, Lemma 3, and Theorem 1 into the Ellentuck structure. The version of Lemma 2 has already been proven in \cite{AFP}, so we present it as Fact 7.

In the proof of Lemma 5, during the construction of trees, at each step, we divide them into two subsets. However, it is possible to divide them into finitely many subsets, with no effect on the outcome. The reason for this change is to make the proof more readable.
	
	\begin{lemma}
		If $\{M_i \colon i \in I\}$ of sets $M_i \subseteq [\omega]^\omega$ is a disjoint $CR-$additive family contained in $\mathbb{NR}$, then $\bigcup\{M_i \colon i \in I\}$  belongs to $\mathbb{NR}$.
	\end{lemma}
	
	\begin{proof}
		Suppose that $\bigcup\{M_i \colon i \in I\}$ does not belong to $\mathbb{NR}$.
		Let $M \subseteq \bigcup\{M_i \colon i \in I\}$ be  a $CR-$set. Let $\{[a_\alpha, A_\alpha] \colon \alpha \in 2^\omega\}$ be a family of all $EL-$sets such that $[a_\alpha, A_\alpha] \subseteq M$.  By Fact 5, there exist distinct $i_\alpha, j_\alpha \in \{i_\beta, j_\beta \in I \colon \beta < \alpha\}$ such that 
		$$[a_\alpha, A_\alpha] \cap[M_{i_\alpha}] \not = \emptyset \textrm{ and }
		[a_\alpha, A_\alpha] \cap[M_{j_\alpha}] \not = \emptyset.$$ 
		Let $J = \{j_\alpha \colon \alpha \in m^\omega\}$. Then, for all $\alpha \in 2^\omega$ $$[a_\alpha, A_\alpha] \cap \bigcup\{M_i \colon i \in J\}  \not = \emptyset $$ and $$[a_\alpha, A_\alpha] \cap  \bigcup\{M_i \colon i \in I \setminus J\}  \not = \emptyset.$$
		Hence $\bigcup\{M_i \colon i \in I\}$ is not a $CR-$set - a contradiction.
	\end{proof}
	\\
	
The next lemma is the "one-dimensional" version of the main result in this subsection given in Theorem 2.	
	
	\begin{lemma}
		Let $f_n \colon [\omega]^\omega \to \{0,1\}$ be a sequence of $CR-$measurable functions. Then there is $[a, A]$ and a subsequence $\{f_{n_k}\}$ which is uniformly convergent on $[a, A]$.
	\end{lemma}
	
	\begin{proof}
		Using Fact 7, it is enough to assume that $\{f_n\}$ are continuous on some $EL-$set $[b, B']$. By Fact 1, one can assume that $\{f_n\}$ are continuous on $[\emptyset, B]$, where $B = b \cup B'$. 
		
		Now, we will construct a family of $EL-$subsets $[a^{\varepsilon_{n_k}}_h, A^{\varepsilon_{n_k}}_h], h \in 2^k, \varepsilon \in \{0,1\}$ of $[\emptyset, B]$ and a subsequence $\{f^{\varepsilon_{n_k}}_{n_k}\}$ of $\{f_n\}$ fulfilling the following properties: for all distinct $h, h' \in 2^k$
		\begin{itemize}
			\item [(1)] $[a^{\varepsilon_{n_k}}_{h}, A^{\varepsilon_{n_k}}_{h}] \subseteq [a^{\varepsilon_{n_{k-1}}}_{h|(k-1)}, A^{\varepsilon_{n_{k-1}}}_{h|(k-1)}]$;
			\item [(2)] $a^{\varepsilon_{n_k}}_h\cap a^{\varepsilon_{n_k}}_{h'} =\emptyset$ and $A^{\varepsilon_{n_k}}_h \cap A^{\varepsilon_{n_k}}_{h'} = \emptyset$; 
			\item [(3)]  $[a^{\varepsilon_{n_k}}_h, A^{\varepsilon_{n_k}}_h] \subseteq (f^{\varepsilon_{n_k}}_{n_k})^{-1}(\{\varepsilon_{n_k}\})$. 
		\end{itemize}
	
	The notation $\varepsilon_{n_k}$ in the conditions above is essential because $\varepsilon_{n_k}$ can be differnt at each level, i.e. $(\varepsilon_{n_k})_{n_k \in \omega}$ form a $0-1$ sequence.
		
		Let $f^{\varepsilon_{n_0}}_{n_0}$ be an arbitrary element of $\{f_n\}$.   The first and inductive steps  are the same. 
		
		Assume that for some $k\in \omega$ we have constructed the subsequence 
		$\{f^{\varepsilon_{n_k}}_{n_k}\}$ of $\{f_n\}$ and the family of $EL-$sets 
		$$\{[a^{\varepsilon_{n_k}}_h, A^{\varepsilon_{n_k}}_h]\colon h \in 2^k, {\varepsilon_{n_k}}\in \{0,1\}\}$$ fulfilling properties $(1) - (3)$. At least one of  $[a^{\varepsilon_{n_k}}_h, A^{\varepsilon_{n_k}}_h]$ is nonempty, hence we continue the construction only for such nonempty subsets.
		
		Fix $h \in 2^k$ such that $[a^{\varepsilon_{n_k}}_h, A^{\varepsilon_{n_k}}_h]$ is nonempty. 
		Choose  $\{f^{\varepsilon_{n_{k+1}}}_{n_{k+1}}\}$ from $\{f_n\}$ such that $$(a^{\varepsilon_{n_k}}_h\cup A^{\varepsilon_{n_k}}_h) \cap (f^{\varepsilon_{n_{k+1}}}_{n_{k+1}})^{-1}(\{\varepsilon_{n_{k}}\}) \not = \emptyset.$$
		Enumerate all subsets of $a^{\varepsilon_{n_k}}_h$ by $a^{\varepsilon_{n_k}}_{h_l}$ for $l=1, 2, ... m$ where $ m = 2^k$. Construct the sequences of subsets of $A^{\varepsilon_{n_k}}_h$:
		$$C^{\varepsilon_{n_k}}_{h_0} \supseteq C^{\varepsilon_{n_k}}_{h_1}\supseteq ... \supseteq C^{\varepsilon_{n_k}}_{h_m}
		\textrm{ and }
		D^{\varepsilon_{n_k}}_{h_0} \supseteq D^{\varepsilon_{n_k}}_{h_1}\supseteq ... \supseteq D^{\varepsilon_{n_k}}_{h_m}$$ 
		as follows: let $C^{\varepsilon_{n_k}}_{h_0},  D^{\varepsilon_{n_k}}_{h_0} \subseteq [A^{\varepsilon_{n_k}}_{h}\setminus a^{\varepsilon_{n_k}}_h]^\omega $ and $C^{\varepsilon_{n_k}}_{h_0}\cap D^{\varepsilon_{n_k}}_{h_0}  = \emptyset$.
		Given $C^{\varepsilon_{n_k}}_{h_l}, D^{\varepsilon_{n_k}}_{h_l}$, whenever there exist $C\subseteq C^{\varepsilon_{n_k}}_{h_l}$ and $D\subseteq D^{\varepsilon_{n_k}}_{h_l}$ such that
		$$[a^{\varepsilon_{n_k}}_{h_l}, C^{\varepsilon_{n_k}}_{h_l}] \subseteq (f^{\varepsilon_{n_{k+1}}}_{n_{k+1}})^{-1}(\{\varepsilon_{n_{k}}\})$$
		and
		$$[a^{\varepsilon_{n_k}}_{h_l}, D^{\varepsilon_{n_k}}_{h_l}] \subseteq (f^{\varepsilon_{n_{k+1}}}_{n_{k+1}})^{-1}(\{\varepsilon_{n_{k}}\}),$$
		then $C^{\varepsilon_{n_k}}_{h_{l+1}}=C^{\varepsilon_{n_k}}_{h_{m}}$ and $D^{\varepsilon_{n_k}}_{h_{l+1}}=D^{\varepsilon_{n_k}}_{h_{m}}$. If not, then then we take $C^{\varepsilon_{n_k}}_{h_{l+1}}=C^{\varepsilon_{n_k}}_{h_l}$ and $D^{\varepsilon_{n_k}}_{h_{l+1}}=D^{\varepsilon_{n_k}}_{h_l}$.
		
		Take
		$$[a^{\varepsilon_{n_{k+1}}}_{h^\smallfrown 0}, A^{\varepsilon_{n_{k+1}}}_{h^\smallfrown 0}] \subseteq [a^{\varepsilon_{n_{k}}}_{h_m} \cup \{\min C^{\varepsilon_{n_{k}}}_{h_m}\}, C^{\varepsilon_{n_{k}}}_{h_m} \setminus \{\min C^{\varepsilon_{n_{k}}}_{h_m}\}],$$  
		$$[a^{\varepsilon_{n_{k+1}}}_{h^\smallfrown 1}, A^{\varepsilon_{n_{k+1}}}_{h^\smallfrown 1}] \subseteq [a^{\varepsilon_{n_{k}}}_{h_m} \cup \{\min D^{\varepsilon_{n_{k}}}_{h_m}\}, D^{\varepsilon_{n_{k}}}_{h_m} \setminus \{\min D^{\varepsilon_{n_{k}}}_{h_m}\}],$$
		$$[a^{\varepsilon_{n_{k+1}}}_{h^\smallfrown 0}, A^{\varepsilon_{n_{k+1}}}_{h^\smallfrown 0}] \subseteq (f^{\varepsilon_{n_{k+1}}}_{n_{k+1}})^{-1}(\{\varepsilon_{n_{k+1}}\})$$ $$[a^{\varepsilon_{n_{k+1}}}_{h^\smallfrown 1}, A^{\varepsilon_{n_k}}_{h^\smallfrown 1}] \subseteq (f^{\varepsilon_{n_{k+1}}}_{n_{k+1}})^{-1}(\{\varepsilon_{n_{k+1}}\}).$$ 
		
		Thus, there are constructed the subsequence $\{f^{\varepsilon_{n_{k}}}_{n_k}\}$ and $$\{[a^{\varepsilon_{n_{k}}}_h, A^{\varepsilon_{n_{k}}}_h]\colon h \in 2^k, {\varepsilon_{n_k}}\in \{0,1\}\}$$ fulfilling $(1)-(3)$. 
		
		Now, take
		a set $$[\emptyset, M^{(\varepsilon_{n_k})_{n_k \in \omega}}] = \bigcap_{k \in \omega} \bigcup_{h \in 2^k} [a^{\varepsilon_{n_k}}_h, A^{\varepsilon_{n_k}}_h]$$
		By Fact 9, the set $[\emptyset ,M^{(\varepsilon_{n_k})_{n_k \in \omega}}]$  is an $EL-$set.
		Let 
		$$[\emptyset, M^\varepsilon] = \bigcup \{[\emptyset, M^{(\varepsilon_{n_k})_{n_k \in \omega}}] \colon  (\varepsilon_{n_k})_{n_k \in \omega} \to \varepsilon\},$$
		where $\varepsilon \in \{0,1\}$.
		One of the  $EL-$sets
		$[\emptyset, M^0]$ or $[\emptyset, M^1]$ is nonempty. Then $$\{f^{\varepsilon_{n_k}}_{n_k}\colon \varepsilon_{n_k} = 0\}$$ is uniformly convergent on $[\emptyset, M^0]$ or $$\{f^{\varepsilon_{n_k}}_{n_k}\colon \varepsilon_{n_k} = 1\}$$ is uniformly convergent on $[\emptyset, M^1]$.
	\end{proof}
	\\
	
	\begin{theorem}
		Let $f_n \colon ([\omega]^\omega)^\omega \to [0,1]$ be a sequence of $CR-$measurable functions. Then, there are $EL-$sets $\{[a_k, A_k]\colon  k \in \omega\}$ and a subsequence $$\{f_{n_k} \colon k \in \omega\}$$ which is uniformly convergent on $\bigotimes_{k < \omega} [a_k, A_k]$.
	\end{theorem}
	
	\begin{proof}
		Since $f_n \colon ([\omega]^\omega)^\omega \to [0,1]$ is $CR$-measurable for all $n \in \omega$, there exist $CR$-measurable functions $g_n \colon K^\omega \to {0,1}^\omega$ such that $f_n = g_n \circ \varphi$, where $\varphi \colon {0,1}^\omega \to [0,1]$ is the canonical continuous "onto" mapping (when one considers elements of $[0,1]$ in their binary expansions).
				Thus, to complete the proof it is enough to  find $EL-$sets  $[a_k, A_k], k \in \omega$ and a subsequence $\{g_{n_k}\}$ which is uniformly convergent on $\bigotimes_{k\in \omega} [a_k, A_k]$.
		We can consider each $g_n$ to be in the form $g_n = (h^k_n \colon k\in \omega)$, where $h^k_n \colon [\omega]^\omega \to \{0,1\}$ are $CR-$measurable. By Fact 7, for each $h^k_n$ and for each $[a, A]$ there is $A'\subseteq A$ such that $A' \subseteq [\omega]^\omega$ and $h^k_n|_{[a, A']}$ is continuous. By Fact 1, we can consider $[\emptyset, B]$ instead of $[a, A']$, where $B=a\cup A'.$ 
		
		Now, consider a partition $\{B_k \colon k \in \omega\}$ of $\omega$ such that $h^k_n|_{[\emptyset, B_k]}$, (or, in the case that it is impossible $h^k_n|_{[\emptyset, B_k]}$ for some infinite $B'_k \subset B_k$). 
		
		We will construct a family $\{\{g_{m,i}\} \colon m \in D_i, i< \omega\}$ (for some set of indices $D_i\subseteq \omega$ of the subsequences $\widehat{h}^i_{m}$ chosen at each step in the construction) of infinite subsequences  of $\{g_n\}$ in the following way.
		Using Lemma 5 to the sequence $\{h^0_n\}$, (i.e. the first coordinates of $\{g_n\}$) we obtain an $EL-$set $[\emptyset, M^{\varepsilon_0}_0]$ and a subsequence $\{\widehat{h}^0_{m}\}$ which is uniformly convergent on $[\emptyset, M^{\varepsilon_0}_0]$, where $\varepsilon_0 \in \{0,1\}$.
		Then, we obtain
		$$g_{m,0} = (\widehat{h}^0_{m}|_{[\emptyset, M^{\varepsilon_0}_0]}, h^1_m, ..., h^k_m,... ).$$
		Now, apply Lemma 5 to the sequence $\{h^1_m\}$, (the second coordinates of $g_{m,0}$) and similar as above we obtain an $EL-$set $[\emptyset, M^{\varepsilon_1}_1]$,  and a subsequence $\{\widehat{h}^1_{m}\}$ which is uniformly convergent on $[\emptyset, M^{\varepsilon_1}_1]$, where $\varepsilon_1 \in \{0,1\}$.
		Thus, we define the subsequence
		$\{g_{m,1}\}$ of $\{g_{m,0}\}$ which is of the form
		$$g_{m,1} = ( \widehat{h}^0_{m}|_{[\emptyset, M^{\varepsilon_0}_0]}, \widehat{h}^1_m|_{[\emptyset, M^{\varepsilon_1}_1]}, h^2_m, ..., h^k_m,... ).$$ 
		Applying Lemma 5 $\omega$ times for each $\{h^k_m\}$ we obtain the family $$\{\{g_{m,i}\} \colon m \in D_i, i\in \omega\}$$ of infinite subsequences of the following properties:
		\begin{itemize}
			\item[(a)] $\{g_{m, i}\}$ is a subsequence of $\{g_{m, i-1}\}$, i.e. $D_i \subseteq D_{i-1}$ for $i \in \omega$,
			\item[(b)] $g_{m,i} = (\bar{h}^k_m \colon k< \omega)$ such that
			$$\bar{h}^k_m = \left\{ \begin{array}{ll}
			\widehat{h}^l_{m}|_{[\emptyset, K^{\varepsilon_l}_l]} & \textrm{for $l\leqslant i$}\\
			h^l_m & \textrm{for $l>i$}\\
			\end{array} \right.,$$
			where $\{\widehat{h}^l_{m}\}$ is the subsequence of $\{\widehat{h}^{l-1}_{m}\}$ uniformly convergent on $[\emptyset, M^{\varepsilon_l}_l],$ $ \varepsilon_l \in \{0,1\}$ and $h^l_m$ is an $l-$th coordinate of $g_{m, l}$.
		\end{itemize}	
		Now, from each subsequence $\{g_{m,i}\}$ one can choose its $i-$th element and define the sequence $\{g_{i,i} \colon i \in \omega\}$. Divide $\{g_{i,i}\}$ into two sequences (one of them may not exist):
		$\{g^0_{i,i}\}$ iff $\{\bar{h}^i_i\}$ are convergent uniformly on $ [\emptyset, M^{\varepsilon_i}_i], \varepsilon_i = 0$ and 	$\{g^1_{i,i}\}$ iff $\{\bar{h}^i_i\}$ are convergent uniformly on  $ [\emptyset, M^{\varepsilon_i}_i], \varepsilon_i = 1$. Then, $\{g^0_{i,i}\}$ is uniformly convergent on $\bigotimes_{i< \omega} [\emptyset, M^0_i]$ or $\{g^1_{i,i}\}$ is uniformly convergent on $\bigotimes_{i< \omega} [\emptyset, M^{1}_i]$. The proof is complete.
	\end{proof}
	\\
	 
	\textbf{Acknowledgments} The author is very grateful to the reviewer for their insightful study of the results presented here and valuable comments and remarks, which allowed the omission of inaccuracies and errors in the paper.

		\noindent
		{\sc Joanna Jureczko}
		\\
		Wroc\l{}aw University of Science and Technology, Poland
		\\ 
		{\sl e-mail: joanna.jureczko@pwr.edu.pl}
	\end{document}